\documentclass[10pt]{amsart}

\usepackage{amsmath, amsfonts, amssymb, amsthm, mathdots, mathrsfs}

\numberwithin{equation}{section}

\begin{document}

\title[$\mathcal{LS}$ method in positive characteristic]{The $\mathcal{LS}$ method for the classical groups \\ in positive characteristic and \\ the Riemann Hypothesis}

\author{Luis Alberto Lomel\'i}

\begin{abstract}
   We provide a definition for an extended system of $\gamma$-factors for products of generic representations $\tau$ and $\pi$ of split classical groups or general linear groups over a non-archimedean local field of characteristic $p$. We prove that our $\gamma$-factors satisfy a list of axioms (under the assumption $p \neq 2$ when both groups are classical groups) and show their uniqueness (in general). This allows us to define extended local $L$-functions and root numbers. We then obtain automorphic $L$-functions $L(s,\tau \times \pi)$, where $\tau$ and $\pi$ are globally generic cuspidal automorphic representations of split classical groups or general linear groups over a global function field. In addition to rationality and the functional equation, we prove that our automorphic $L$-functions satisfy the Riemann Hypothesis.
\end{abstract}

\maketitle

\section*{introduction}

Let ${\bf G}_m$ and ${\bf G}_n$ denote either split classical groups or general linear groups of ranks $m$ and $n$, respectively. Let $k$ be a global function field with finite field of constants $\mathbb{F}_q$ and ring of ad\`eles $\mathbb{A}_k$. We present a theory of automorphic $L$-functions $L(s,\tau \times \pi)$, where $\tau$ and $\pi$ are globally generic cuspidal automorphic representations of ${\bf G}_m(\mathbb{A}_k)$ and ${\bf G}_n(\mathbb{A}_k)$.

The case of ${\bf G}_m = {\rm GL}_m$ and ${\bf G}_n$ a classical group is made possible by our work on the Langlands-Shahidi method in positive characteristic for the classical groups \cite{lomeli2009, lomelipreprint}. Already in this case, the $\mathcal{LS}$ method has particularly interesting applications. In addition to the Ramanujan conjecture for the classical groups over function fields established in \cite{lomeli2009}, the zeros of $L(s,\tau \times \pi)$ satisfy $\Re(s) = 1/2$.

We note that the case of ${\bf G}_m = {\rm GL}_m$ and ${\bf G}_n = {\rm GL}_n$ gives rise to Rankin-Selberg factors. Indeed, we include a treatise of this case in a self contained manner within the Langlands-Shahidi method in \cite{lomelipreprint}. And, we provide a short proof of the equality of local factors when defined via different methods in \cite{hlappendix}. Thanks to the work of Lafforgue on the Langlands correspondence for ${\rm GL}_N$ over function fields, the Riemann Hypothesis is available for $L$-functions of products of cuspidal automorphic representations of ${\rm GL}_m(\mathbb{A}_k)$ and ${\rm GL}_n(\mathbb{A}_k)$ \cite{laff2002}. In this very important case, all of our results are available with no restriction on the characteristic of $k$.

Notice that the case of two classical groups ${\bf G}_m$ and ${\bf G}_n$ is treated for the first time here in positive characteristic. We provide axioms for an extended system of $\gamma$-factors, $L$-functions and root numbers which cover all of the above mentioned cases for ${\bf G}_m$ and ${\bf G}_n$. We first establish existence and uniqueness of $\gamma$-factors, Theorem~1.5, and then include existence and uniqueness of local $L$-functions and $\varepsilon$-factors in Theorem~4.3. The theory is complete under the assumption ${\rm char}(k) \neq 2$.

We begin by introducing notation that is useful when dealing with systems of $\gamma$-factors, $L$-functions and root numbers. Local factors $\gamma(s,\tau \times \pi,\psi)$, $L(s,\tau \times \pi)$ and $\varepsilon(s,\tau \times \pi,\psi)$ are defined on the local class $\mathfrak{ls}(p,{\bf G}_m,{\bf G}_n)$, while global $L$-functions and root numbers
\begin{equation*}
   L(s,\tau \times \pi)  = \prod_v L(s,\tau_v \times \pi_v), \quad \varepsilon(s,\tau \times \pi) = \prod_v \varepsilon(s,\tau_v \times \pi_v,\psi_v)
\end{equation*}
are defined on the global class $\mathcal{LS}(p,{\bf G}_m,{\bf G}_n)$ (see \S~1.1 and \S~1.2). We often write $\mathfrak{ls}(p)$ and $\mathcal{LS}(p)$ when there is no need to specify the groups ${\bf G}_m$ and ${\bf G}_n$. Theorem~4.4 can then be succinctly stated as follows:

\subsection*{Theorem} \emph{Automorphic $L$-functions on $\mathcal{LS}(p)$, $p \neq 2$, satisfy the following properties:}
\begin{itemize}
   \item[(i)] (Rationality). \emph{$L(s,\tau \times \pi)$ converges on a right half plane and has a meromorphic continuation to a rational function on $q^{-s}$.}
   \item[(ii)] (Functional equation). $L(s,\tau \times \pi) = \varepsilon(s,\tau \times \pi) L(1-s,\tilde{\tau} \times \tilde{\pi})$.
   \item[(iii)] (Riemann Hypothesis) \emph{The zeros of $L(s,\tau \times \pi)$ are contained in the line $\Re(s) = 1/2$}.
\end{itemize}

\medskip

Let us now give a more detailed description of the contents of the article. Theorem~1.5 concerns the existence an uniqueness of a system of $\gamma$-factors on $\mathfrak{ls}(p)$. The theorem is true with no assumption on $p$ for $\mathfrak{ls}(p,{\rm GL}_m,{\bf G}_n)$. However, we assume $p \neq 2$ in order to produce $\gamma$-factors on $\mathfrak{ls}(p,{\bf G}_m,{\bf G}_n)$ when both ${\bf G}_m$ and ${\bf G}_n$ are classical groups. Extended $\gamma$-factors satisfy several local properties including a twin multiplicativity property when the representations are obtained via parabolic induction, as well as a stability property on $\mathfrak{ls}(p,{\rm GL}_1,{\bf G}_n)$ with respect to highly ramified characters of ${\rm GL}_1$. Globally, $\gamma$-factors make an important appearance in the functional equation for partial $L$-functions on $\mathcal{LS}(p)$.

All of \S~2 is devoted to a proof of the existence part of Theorem~1.5. With the Langlands-Shahidi method now complete for the split classical groups in positive characteristic, we give a straightforward and linearly ordered presentation of $\gamma$-factors on $\mathfrak{ls}(p,{\rm GL}_m,{\bf G}_n)$. We recall the basic definitions in \S~2.1. Then, our results on the Siegel Levi case for the split classical groups \cite{lomelipreprint} are summarized in \S~2.2. The Siegel Levi case allows us to define exterior and symmetric square $\gamma$-factors, which are uniquely characterized and are proved to be in accordance with the local Langlands conjecture for ${\rm GL}_m$ in \cite{hl2011}. Filling any gaps that were left in \cite{lomeli2009}, the case of $\mathfrak{ls}(p,{\rm GL}_m,{\bf G}_n)$ is presented in \S~2.3 with no assumption on the characteristic. The new case of $\mathfrak{ls}(p,{\bf G}_m,{\bf G}_n)$, with both ${\bf G}_m$ and ${\bf G}_n$ classical groups, is developed in \S~2.4 under the assumption $p \neq 2$.

There are fundamental differences in local-to-global arguments between global fields of characteristic zero and characteristic $p$. In \S~3, we use a variation of a result of Vign\'eras that allows us to lift a local pair of irreducible supercuspidal generic representations $\tau_0$ and $\pi_0$ to a global pair of globally generic cuspidal automorphic representations $\tau$ and $\pi$ (see Proposition~3.1). We remark that, over number fields, Shahidi makes a crucial improvement upon the local-to-global argument of Henniart and Vign\'eras by incorporating the archimedean theory that is available at infinite places (Proposition~5.1 of \cite{sha1990}). Over function fields, the main difference is due to the fact that all places are non-archimedean; a place at infinity plays the role of the archimedean places. As a further remark, in the case of two general linear groups we have a much more precise local-to-global result \cite{hl2011}, which allows us to remove stability from the list of properties required in the characterization \cite{hlappendix}.

The uniqueness part of Theorem~1.5 is proved in \S~3 over a global function field with no restriction on $p$. We first treat the case of $\mathfrak{ls}(p,{\rm GL}_1,{\bf G}_n)$ in \S~3.2, where we can use stability of $\gamma$-factors combined with the Grundwald-Wang theorem of class field theory. We then proceed to the general case in \S~3.3. Our method of proof resembles the approach taken in \cite{lomeli2009}, and we refer to the introduction for further remarks on the local-to-global argument over a global function field (see also \S~3).

Building upon extended $\gamma$-factors, we axiomatize local $L$-functions and root numbers on $\mathfrak{ls}(p)$ in \S~4.1. Additional properties of $\gamma$-factors are recorded in \S~4.2, these include a local functional equation for which we provide a proof. We extend Theorem~1.5 to a theorem on extended $\gamma$-factors, local $L$-functions and root numbers in \S~4.3. Finally, we establish our main global results in \S~4.4, under the assumption $p \neq 2$. Theorem~4.4 includes rationality, the functional equation and the Riemann Hypothesis for automorphic $L$-functions on $\mathcal{LS}(p)$.

Our general results are possible since we have established a Langlands functorial lift from globally generic cuspidal automorphic representations $\pi$ of a classical group ${\bf G}_n$ to automorphic representations $\Pi$ of ${\rm GL}_N$ \cite{lomeli2009}. The integer $N$ is obtained from the table of \S~1.3, by minimally embedding the connected component of the Langlands dual group ${}^L{G}_n$ of ${\bf G}_n$ into ${\rm GL}_N(\mathbb{C})$. The image of functoriality can be further expressed as an isobaric sum
\begin{equation*}
   \Pi = \Pi_1 \boxplus \cdots \boxplus \Pi_d,
\end{equation*}
where each $\Pi_i$ is a self-dual cuspidal automorphic representation of ${\rm GL}_{n_i}$. They satisfy the additional condition that, for each $i$, a partial $L$-function $L^S(s,\Pi_i,r \circ \rho_{n_i})$ has a pole at $s=1$, where $r = \wedge^2$ or ${\rm Sym}^2$ depending on the classical group \cite{soudry2005}. Exterior and symmetric square $L$-functions, and related $\gamma$-factors, are thoroughly studied in \cite{hl2011,lomelipreprint}. Furhtermore, \cite{lomelipreprint} also develops the necessary theory for the Siegel Levi case of quasi-split Unitary groups. This leads us to Asai $\gamma$-factors and $L$-functions, which are uniquely characterized in \cite{hlpreprint}. They play a similar role for the quasi-split unitary groups, as the exterior and symmetric square $L$-functions do for the split classical groups, when describing the image of the functorial lift of \cite{lomeliunitary} as an isobaric sum. 

I would like to thank Freydoon Shahidi for many insightful mathematical conversations. Thanks are also due to Guy Henniart; our collaborative work greatly influenced the axiomatization of the system of $\gamma$-factors, $L$-functions and $\varepsilon$-factors presented in this paper. I would also like to thank the D\'epartement de Math\'ematiques d'Orsay, Universit\'e Paris-Sud~11, for their hospitality during a short visit to conduct research in May-June of 2012. Part~(iii) of Theorem~4.4 for the case of $\mathcal{LS}(p,{\rm GL}_m,{\bf G}_n)$, ${\bf G}_n$ a classical group, was an unpublished result of the author since the summer of 2009. The opportunity presented itself to write its proof together with the extended case of $\mathcal{LS}(p,{\bf G}_m,{\bf G}_n)$ during this visit.

\section{Extended $\gamma$-factors}

Let ${\bf G}_n$ be either the general linear group ${\rm GL}_n$ or a split classical group of rank $n$. Given a ring $R$ and an algebraic group ${\bf G}$ defined over $R$, we often let $G$ denote its group of rational points. Given a non-archimedean local field $F$, let $\mathcal{O}_F$ denote its ring of integers, $\mathfrak{p}_F$ its maximal ideal, $\varpi_F$ a uniformizer, and $q_F$ the cardinality of its residue field. Given a global function field field $k$, we let $q$ denote the cardinality of its field of constants; for a place $v$ of $k$, we let $q_v$ be the cardinality of the residue field of $k_v$. Given a representation $\sigma$, we let $\tilde{\sigma}$ denote its contragredient.

\subsection{Local notation} Let $\mathfrak{ls}(p,{\bf G}_m,{\bf G}_n)$ denote the class of quadruples $(F,\tau,\pi,\psi)$ consisting of: a non-archimedean local field $F$ of characteristic $p$; irreducible generic representations $\tau$ of $G_m = {\bf G}_m(F)$ and $\pi$ of $G_n = {\bf G}_n(F)$; and, a non-trivial character $\psi$ of $F$.

Given $(F,\tau,\pi,\psi) \in \mathfrak{ls}(p,{\bf G}_m,{\bf G}_n)$, we call it tempered (resp. discrete series, supercuspidal) if $\tau$ and $\pi$ are tempered representations (resp. discrete series, supercuspidal).

\subsection{Global notation} Let $\mathcal{LS}(p,{\bf G}_m,{\bf G}_n)$ denote the class of quadruples $(k,\tau,\pi,\psi)$ consisting of: a global function field $k$ of characteristic $p$; globally generic cuspidal automorphic representations $\tau = \otimes' \, \tau_v$ of $G_m = {\bf G}_m(\mathbb{A}_k)$ and $\pi = \otimes' \, \pi_v$ of $G_n = {\bf G}_n(\mathbb{A}_k)$; and, a non-trivial character $\psi = \otimes \,\psi_v$ of $k \backslash \mathbb{A}_k$.

\subsection*{Remark} We often write $\mathfrak{ls}(p)$ and $\mathcal{LS}(p)$ when there is no need to specify the groups ${\bf G}_m$ and ${\bf G}_n$.

\subsection{$L$-groups, principal series and partial $L$-functions} The connected components of the $L$-groups for the split classical groups are embedded into an appropriate dual group of ${\rm GL}_N$ according to the following table:

\medskip

\begin{center}
\begin{tabular}{|c|c|c|} \hline
   ${\bf G}_n$			& ${}^LG_n^0 \hookrightarrow {}^L{\rm GL}_N^0$							& ${\rm GL}_N$ \\ \hline
   ${\rm SO}_{2n+1}$	& ${\rm Sp}_{2n}(\mathbb{C}) \hookrightarrow {\rm GL}_{2n}(\mathbb{C})$		& ${\rm GL}_{2n}$ \\
   ${\rm SO}_{2n}$		& ${\rm SO}_{2n}(\mathbb{C}) \hookrightarrow {\rm GL}_{2n}(\mathbb{C})$		&${\rm GL}_{2n}$ \\
   ${\rm  Sp}_{2n}$		& ${\rm SO}_{2n+1}(\mathbb{C}) \hookrightarrow {\rm GL}_{2n+1}(\mathbb{C})$	& ${\rm GL}_{2n+1}$ \\ \hline
\end{tabular}
\end{center}

\medskip

We include the possibility of ${\bf G}_n = {\rm GL}_n$, where we take $N = n$. Also, we let $\rho_n$ denote the standard representation of ${\rm GL}_n(\mathbb{C})$.

Let $(F,\tau,\pi,\psi) \in \mathfrak{ls}(p,{\bf G}_m,{\bf G}_n)$ and assume that $\tau$ and $\pi$ are unramified principal series. Then, the Satake parametrization gives semisimple conjugacy classes $\left\{ A_\tau \right\}$ in ${}^LG_m^0 \hookrightarrow {\rm GL}_M(\mathbb{C})$ and $\left\{ B_\pi \right\}$ in ${}^LG_n^0 \hookrightarrow {\rm GL}_N(\mathbb{C})$. Then, $L$-functions for unramified principal series representations are defined by
\begin{equation*}
   L(s,\tau \times \pi) = \dfrac{1}{\det(I - \rho_M(A_\tau) \otimes \rho_N(B_\pi)q_F^{-s})}.
\end{equation*}

Given $(k,\tau,\pi,\psi) \in \mathcal{LS}(p,{\bf G}_m,{\bf G}_n)$, we take $S$ to be a finite set of places of $k$ such that $\tau$, $\pi$ and $\psi$ are unramified for $v \notin S$. The corresponding partial $L$-functions are defined by
   \begin{equation*}
      L^S(s,\tau \times \pi) = \prod_{v \notin S} L(s,\tau_v \times \pi_v).
   \end{equation*}
We can show that partial $L$-functions converge on a right half plane; in fact, they have a meromorphic continuation to a rational function on $q^{-s}$.

\subsection{Axioms for a system of $\gamma$-factors}

The Langlands-Shahidi method in positive characteristic allows us to produce a system of rational functions $\gamma(s,\tau \times \pi,\psi) \in \mathbb{C}(q_F^{-s})$ on $\mathfrak{ls}(p,{\rm GL}_m,{\bf G}_n)$. In this article, we concoct a system of extended $\gamma$-factors on $\mathfrak{ls}(p,{\bf G}_m,{\bf G}_n)$. Extended $\gamma$-factors can be characterized by a list of local properties together with their role in the global functional equation.

\begin{itemize}
   \item[(i)] (Naturality). \emph{Let $(F,\tau,\pi,\psi) \in \mathfrak{ls}(p)$, and let $\eta: F' \rightarrow F$ be an isomorphism of local fields. Let $(F',\tau',\pi',\psi') \in \mathfrak{ls}(p)$ be the quadruple obtained from $(F,\tau,\pi,\psi)$ via $\eta$. Then
   \begin{equation*}
      \gamma(s,\tau \times \pi,\psi) = \gamma(s,\tau' \times \pi',\psi').
   \end{equation*}
   }
   \item[(ii)] (Isomorphism). \emph{Let $(F,\tau,\pi,\psi) \in \mathfrak{ls}(p)$. If $(F,\tau',\pi',\psi) \in \mathfrak{ls}(p)$ is such that $\tau \simeq \tau'$ and $\pi \simeq \pi'$, then
   \begin{equation*}
      \gamma(s,\tau \times \pi,\psi) = \gamma(s,\tau' \times \pi',\psi).
   \end{equation*}
   }
   \item[(iii)] (Compatibility with class field theory). \emph{Let $(F,\chi_1,\chi_2,\psi) \in \mathfrak{ls}(p,{\rm GL}_1,{\rm GL}_1)$. Then
   \begin{equation*}
      \gamma(s,\chi_1 \times \chi_2,\psi) = \gamma(s,\chi_1 \chi_2,\psi).
   \end{equation*}
   }
   \item[(iv)] (Multiplicativity). \emph{Let $(F,\tau_i,\pi_j,\psi) \in \mathfrak{ls}(p,{\bf G}_{m_i},{\bf G}_{n_j})$, $0 \leq i \leq d$, $0 \leq j \leq e$; ${\bf G}_{m_0}$ and ${\bf G}_{n_0}$ can be classical groups or general linear groups; ${\bf G}_{m_i} = {\rm GL}_{m_i}$ and ${\bf G}_{n_j} = {\rm GL}_{n_j}$ for $1 \leq i \leq d$, $1 \leq j \leq e$. Set $m = \sum m_i$ and $n = \sum n_j$. If ${\bf G}_{m_0}$ (resp. ${\bf G}_{n_0}$) is a classical group, take ${\bf G}_m$ (resp. ${\bf G}_n$) to be a classical group of the same type. Let ${\bf P}_m$ (resp. ${\bf P}_n$) be the standard parabolic subgroup of ${\bf G}_m$ (resp. ${\bf G}_n$) with Levi ${\bf M}_m = \prod_{i=1}^{d} {\rm GL}_{m_i} \times {\bf G}_{m_0}$ (resp. ${\bf M}_n = \prod_{j=1}^{e} {\rm GL}_{n_j} \times {\bf G}_{n_0}$). First, assume $m_0 \geq 1$ and $n_0 \geq 1$. Let $\tau$ be the generic constituent of
   \begin{equation*}
      {\rm ind}_{P_m}^{G_m}(\tau_1 \otimes \cdots \otimes \tau_d \otimes \tau_0),
   \end{equation*}
and let $\pi$ be the generic constituent of
   \begin{equation*}
      {\rm ind}_{P_n}^{G_n}(\pi_1 \otimes \cdots \otimes \pi_e \otimes \pi_0).
   \end{equation*}
When $m_0=0$ (resp. $n_0=0)$ we make the following conventions: take $\tau_0$ (resp. $\pi_0$) to be the trivial character of ${\rm GL}_1(F)$ if ${\bf G}_m$ (resp. ${\bf G}_n$) is a symplectic group; in all other cases, we interpret $\gamma$-factors involving $\tau_0$ (resp. $\pi_0$) to be trivial.
   \begin{itemize}
   \item[(iv.a)] If both ${\bf G}_m$ and ${\bf G}_n$ are classical groups, then
   \begin{align*}
      \gamma(s,\tau \times \pi,\psi) &= \gamma(s,\tau_0 \times \pi_0) \\
      	  & \times \prod_{i=1}^d \gamma(s,\tau_i \times \pi_0,\psi) \gamma(s,\tilde{\tau}_i \times \pi_0,\psi)
      	     \prod_{j=1}^e \gamma(s,\tau_0 \times \pi_j,\psi) \gamma(s,\tau_0 \times \tilde{\pi}_j,\psi) \\
	  &    \times \prod_{1 \leq h \leq d, 1 \leq l \leq e} \gamma(s,\tau_h \times \pi_l,\psi) \gamma(s,\tau_h \times \tilde{\pi}_l,\psi) \gamma(s,\tilde{\tau}_h \times \pi_l,\psi) \gamma(s,\tilde{\tau}_h \times \tilde{\pi}_l,\psi).
   \end{align*}
   \item[(iv.b)] If ${\bf G}_m = {\rm GL}_m$ and ${\bf G}_n$ is a classical group, then
   \begin{equation*}
      \gamma(s,\tau \times \pi,\psi) = \prod_{i=0}^d \gamma(s,\tau_i \times \pi_0,\psi) 
      							\times \prod_{i=1}^d \prod_{j=1}^e \gamma(s,\tau_i \times \pi_j,\psi) \gamma(s,\tau_i \times \tilde{\pi}_j,\psi).
   \end{equation*}
   \item[(iv.c)] If ${\bf G}_m = {\rm GL}_m$ and ${\bf G}_n = {\rm GL}_n$, then
   \begin{equation*}
      \gamma(s,\tau \times \pi,\psi) = \prod_{i,j} \gamma(s,\tau_i \times \pi_j,\psi).
   \end{equation*}
   \end{itemize}
   }
   \item[(v)] (Dependence on $\psi$). \emph{Let $(F,\tau,\pi,\psi) \in \mathfrak{ls}(p,{\bf G}_m,{\bf G}_n)$. Given $a \in F^\times$, let $\psi^a$ denote the character of $F$ given by $\psi^a(x) = \psi(ax)$. Then
   \begin{equation*}
      \gamma(s,\tau \times \pi,\psi^a) = \omega_\tau(a)^h \omega_\pi(a)^l \left| a \right|_F^{hl(s - \frac{1}{2})} \gamma(s,\tau \times \pi,\psi),
   \end {equation*}
   where $h = 2m$ if ${\bf G}_m = {\rm SO}_{2m}$, ${\rm SO}_{2m+1}$; $h = 2m+1$ if ${\bf G}_m = {\rm Sp}_{2m}$; $h = m$ if ${\bf G}_m = {\rm GL}_m$; and, similarly for $l$, depending on ${\bf G}_n$.
   }
   \item[(vi)] (Stability). \emph{Let $(F,\eta,\pi_i,\psi) \in \mathfrak{ls}(p,{\rm GL}_1,{\bf G}_n)$, $i=1$, $2$, where $\pi_1$ and $\pi_2$ have the same central character and $\eta$ is highly ramified. Then
   \begin{equation*}
      \gamma(s,\eta \times \pi_1,\psi) = \gamma(s,\eta \times \pi_2,\psi).
   \end{equation*}
   }
   \item[(vii)] (Functional equation). \emph{Let $(k,\tau,\pi,\psi) \in \mathcal{LS}(p)$, then
   \begin{equation*}
      L^S(s,\tau \times \pi) = \prod_{v \in S} \gamma(s,\tau_v \times \pi_v,\psi_v) L^S(1-s,\tilde{\tau} \times \tilde{\pi}).
   \end{equation*}
   }
\end{itemize}

\subsection{Theorem} \label{gammathm}\emph{There {\bf exists} a system of $\gamma$-factors on $\mathfrak{ls}(p,{\rm GL}_m,{\bf G}_n)$ satisfying properties $(i)-(vii)$. If $p \neq 2$, there {\bf exists} a system of $\gamma$-factors on $\mathfrak{ls}(p,{\bf G}_m,{\bf G}_n)$ satisfying properties $(i)-(vii)$. Any system of $\gamma$-factors satisfying properties $(i)-(vii)$ is {\bf uniquely} determined.}

\section{Existence}

A treatise of $\gamma$-factors, $L$-functions and root numbers for general linear groups is presented in a self contained manner within the Langlands-Shahidi method in \cite{lomelipreprint} and the appendix \cite{hlappendix}. We now complete the study begun in \cite{lomeli2009} for the cases involving split classical groups.

\subsection{The Langlands-Shahidi local coefficient for the split classical groups in positive characteristic} Let ${\bf G}$ be a split classical group of rank $l$ and let ${\bf B} = {\bf T}{\bf U}$ be the Borel subgroup of upper triangular matrices with maximal torus ${\bf T}$ and unipotent radical ${\bf U}$. Let ${\bf P}$ be the standard maximal parabolic subgroup of ${\bf G}$ with maximal Levi ${\bf M} = {\rm GL}_m \times {\bf G}_n$, $l=m+n$, and unipotent radical ${\bf N}$. Given $(F,\tau,\pi,\psi) \in \mathfrak{ls}(p,{\rm GL}_m,{\bf G}_n)$, we can form a smooth irreducible generic representation $\sigma = \tau \otimes \tilde{\pi}$ of $M$.

Let $\Sigma$ denote the roots of ${\bf G}$ with respect to the maximal torus ${\bf T}$, $\Sigma^+$ the positive roots, and $\Delta$ a basis fixed by our choice of Borel subgroup. Let ${\bf N}_\alpha$ denote the one parameter subgroup associated to $\alpha \in \Sigma$. The surjection ${\bf N} \twoheadrightarrow {\bf N} / \prod_{\alpha \in \Sigma^+ - \Delta} {\bf N}_\alpha \simeq \prod_{\alpha \in \Delta} {\bf N}_\alpha$ allows us to extend the character $\psi$ of $F$ to a non-degenerate character $\psi$ of $N$ by setting $\psi(\sum_{\alpha \in \Delta} x_\alpha) = \prod_{\alpha \in \Delta} \psi(x_\alpha)$. Let $\alpha_m \in \Delta$ be such that ${\bf P}$ is the maximal parabolic subgroup corresponding to $\Delta - \left\{ \alpha_m \right\}$. Then, we let $\psi_M$ be the character of $U_M = U \cap M$ obtained from $\psi$ via the surjection ${\bf U}_M \twoheadrightarrow \prod_{\alpha \in \Delta - \left\{ \alpha_m \right\}} {\bf N}_\alpha$.

For every $\alpha \in \Delta$, we fix an embedding for the corresponding semisimple rank one group ${\bf G}_\alpha$ into ${\bf G}$ and fix a representative $w_\alpha$ for the corresponding Weyl group element as in \cite{lomelipreprint}. We abuse notation and identify Weyl group elements with their fixed representatives. Let $w_0 = w_l w_{l,M}$, where $w_l$ is the longest Weyl group element of ${\bf G}$ and $w_{l,M}$ is the longest Weyl group element with respect to ${\bf M}$. Then, the non-degenerate characters $\psi$ of $N$ and $\psi_M$ of $U_M$ are $w_0$-compatible, i.e., $\psi(u) = \psi_M(w_0^{-1}uw_0)$ for $u \in U_M$.

Let $\sigma = \tau \otimes \tilde{\pi}$ be a $\psi_M$-generic representation of $M = {\rm GL}_m(F) \times G_n$. We then let ${\rm I}(s,\sigma)$ be the unitarily induced representation
\begin{equation*}
   {\rm ind}_{P}^{G_l} (\left| {\rm det}(\cdot) \right|_F^{s} \tau \otimes \tilde{\pi}).
\end{equation*}
Let ${\rm V}(s,\sigma)$ denote the space of ${\rm I}(s,\sigma)$. If $\lambda_M$ is a Whittaker functional for ${\sigma}$, then ${\rm I}(s,\sigma)$ is $\psi$-generic for the Whittaker functional $\lambda_\psi$ given by
\begin{equation*}
   \lambda_\psi(s,\sigma)f = \int_N \lambda_M(w_0^{-1}n) \overline{\psi}(n) \, dn,
\end{equation*}
where $f \in {\rm V}(s,\sigma)$. The integral on the right hand side converges as a principal value integral over compact open subgroups of $N$.

We also have an intertwining operator ${\rm A}(s,\sigma,w_0): {\rm V}(s,\sigma) \rightarrow {\rm V}(-s,w_0^{-1}(\sigma))$, given by
\begin{equation*}
   {\rm A}(s,\sigma,w_0)f(g) = \int_N f(w_0^{-1}ng) \, dn,
\end{equation*}
where we write $w_0(\sigma)$ for the representation given by $w_0(\sigma)(x) = \sigma(w_0^{-1}xw_0)$. It converges for $\Re(s) \gg 0$ and extends to a rational operator on $q_F^{-s}$.

The local coefficient is then defined using the uniqueness property of Whittaker models and the relationship
\begin{equation*}
   \lambda_\psi(s,\sigma)f = C_\psi(s,\sigma,w_0) \lambda_\psi(-s,w_0(\sigma)) {\rm A}(s,\sigma,w_0)f.
\end{equation*}
The local coefficient $C_\psi(s,\sigma,w_0)$ is a rational function on $q_F^{-s}$.

\subsection{The case of a Siegel Levi subgroup} The case of a Siegel Levi subgroup ${\bf M} \simeq {\rm GL}_n$ of ${\bf G}_n$ is studied in \cite{hl2011} when ${\bf G}_n = {\rm SO}_{2n}$ or ${\rm SO}_{2n+1}$. In these cases the Langlands-Shahidi method allows us to study exterior square and symmetric square $L$-functions; the case of a Siegel Levi subgroup of ${\rm Sp}_{2n}$ is included in \cite{lomelipreprint}. These cases provide an important step in the Langlands-Shahidi method for the split classical groups. Given an irreducible generic representation $\tau$ of ${\rm GL}_n(F)$, we define
\begin{equation*}
   C_\psi(s,\tau,w_0) = \left\{ \begin{array}{ll}
   						\gamma(s,\tau,{\rm Sym}^2 \rho_n,\psi) & \text{ if } {\bf G} = {\rm SO}_{2n+1} \\
						\gamma(s,\tau, \wedge^2 \rho_n,\psi) & \text{ if } {\bf G} = {\rm SO}_{2n}
				\end{array} . \right.
\end{equation*}

Here, $\rho_n$ is the standard representation of ${\rm GL}_n(\mathbb{C})$, the dual group of ${\rm GL}_n$, and $\gamma(s,\tau,\psi)$ is a Godement-Jacquet $\gamma$-factor. For unramified principal series representations, the above definition agrees with the Satake parametrization and provides a definition of exterior square and symmetric square local factors for general smooth representations. We now state the main result of \cite{hl2011}, which shows that Langlands-Shahidi $\gamma$-factors are in accordance with the local Langlands correspondence in positive characteristic \cite{lrs1993}.

\medskip

\noindent{\bf Theorem.} \emph{Let $\tau$ be an irreducible smooth representation of ${\rm GL}_n(F)$ and let $\sigma$ be an $n$-dimensional ${\rm Frob}$-semisimple $\ell$-adic representation of the Weil-Deligne group in the isomorphism class $\sigma(\tau)$ corresponding to $\tau$ via the local Langlands correspondence. Then
\begin{equation*}
   \gamma(s,\tau,r \circ \rho_n,\psi) = \gamma(s,r \circ \sigma,\psi),
\end{equation*}
where $r = {\rm Sym}^2$ or $\wedge^2$.}

For the remaining case of a Siegel Levi subgroup of a split classical group, let $\gamma(s,\tau,\psi)$ denote a Godement-Jacquet $\gamma$-factor. Then
\begin{equation*}
   C_\psi(s,\tau,w_0) = \gamma(s,\tau,\psi) \gamma(2s,\tau,\wedge^2 \rho_n,\psi) \text{ if } {\bf G} = {\rm Sp}_{2n}.
\end{equation*}

\subsection{The case of $\mathfrak{ls}(p,{\rm GL}_m,{\bf G}_n)$} The study of $\gamma$-factors, $L$-functions and root numbers on $\mathfrak{ls}(p,{\rm GL}_m,{\bf G}_n)$ and $\mathcal{LS}(p,{\rm GL}_m,{\bf G}_n)$ was begun in \cite{lomeli2009}. We now gather the necessary results that establish the existence part of Theorem~\ref{gammathm} in these cases; we use the conventions of \cite{lomelipreprint} regarding Weyl group element representatives and normalization of Haar measures. Let $(F,\tau,\pi,\psi) \in \mathfrak{ls}(p,{\rm GL}_m,{\bf G}_n)$ and let $\sigma = \tau \otimes \tilde{\pi}$. We first assume that $\sigma$ is a $\psi_M$-generic representation of $M$.

Having defined exterior and symmetric square $\gamma$-factors, which are in accordance with the local Langlands correspondence for ${\rm GL}_n$, the $\gamma$-factors $\gamma(s,\tau \times \pi,\psi)$ are defined via the local coefficient:
\begin{equation*}
   C_\psi(s,\sigma,w_0) = \left\{ \begin{array}{ll}
   						\gamma(s,\tau \times \pi,\psi) 
						\gamma(2s,\tau,{\rm Sym}^2 \rho_n,\psi)	& \text{ if } {\bf G} = {\rm SO}_{2n+1} \\
						\gamma(s,\tau \times \pi,\psi)
						\gamma(2s,\tau,\wedge^2 \rho_n,\psi) 	& \text{ if } {\bf G} = {\rm Sp}_{2n} \text{ or } {\rm SO}_{2n}
				\end{array} . \right.
\end{equation*}

Let $(F,\tau,\pi,\psi) \in \mathfrak{ls}(p,{\rm GL}_m,{\bf G}_n)$, with $\sigma = \tau \otimes \tilde{\pi}$ $\psi_M$-generic. An isomorphism of local fields $\eta: F' \rightarrow F$, takes normalized Haar measures on ${\bf N}(F)$ to normalized Haar measures on ${\bf N}(F')$. Hence, the local coefficients $C_\psi(s,\sigma,w_0)$ and $C_{\psi'}(s,\sigma',w_0)$ are equal. Also, we have an isomorphism property for the local coefficient given two isomorphic $\psi_M$-generic representations $\sigma$ and $\sigma'$ of $M$. Thus, properties (i) and (ii) follow for $\psi_M$-generic $\sigma$ and $\sigma'$. Property~(iii) is included in the list of semisimple rank one cases of \cite{lomelipreprint}. Property~(iv.b) for the classical groups is Theorem~6.7 of \cite{lomeli2009}, where it is explicitly stated.

We now discuss the relationship for $\gamma$-factors as the character varies. Given a triple $(F,\tau,\pi,\psi) \in \mathfrak{ls}(p,{\rm GL}_m,{\bf G}_n)$, the representation $\sigma = \tau \otimes \tilde{\pi}$ is generic with respect to a non-degenerate character $\chi_M$ of $U_M$. We will take $\chi$ to be a non-degenerate character of $U$, which is $w_0$-compatible with $\chi_M$. Now, given the group ${\bf G}$, we can embed it into a group $\widetilde{\bf G}$ with Borel subgroup $\widetilde{\bf B} = \widetilde{\bf T}{\bf U}$. The group $\widetilde{\bf G}$ has the same derived group as ${\bf G}$ and has ${\rm H}^1({\rm Z}_{\widetilde{G}}) = \left\{ 1 \right\}$ (see \S~5 of \cite{sha2002}). Thus, $\widetilde{T}$ acts transitively on the set of non-degenerate characters of $U$. Let $t \in \widetilde{T}$ be such that the non-degenerate character $\psi_M$ of $U_M$ is equal to $\chi_{t,M}$, where $\chi_{t,M}(u) = \chi_M(t^{-1}ut)$. The character $\chi$ is taken so that $\psi = \chi_t$ on $U$, and $w_0$-compatibility is preserved for the action of $t \in \widetilde{T}$ on the non-degenerate characters. Let $\sigma_t$ be defined by $\sigma_t(m) = \sigma(t^{-1} m t)$. It is a $\psi_M$-generic representation. Then, we have $\gamma$-factors defined on all of $\mathfrak{ls}(p,{\rm GL}_m,{\bf G}_n)$ via the local coefficient as follows
\begin{equation*}
   \gamma(s,\pi,\psi) = C_\psi(s,\pi_t,w_0).
\end{equation*}
We have a formula for the local coefficient when the character $\psi$ varies. Written explicitly for $\gamma$-factors gives property~(v) on $\mathfrak{ls}(p,{\rm GL}_m,{\bf G}_n)$.

Property~(vi), stability of $\gamma$-factors in positive characteristic, is the content of Theorem~6.12 of \cite{lomeli2009} for the split classical groups. This includes the case of characteristic $2$.

Finally, let $(F,\tau,\pi,\psi) \in \mathcal{LS}(p,{\rm GL}_m,{\bf G}_n)$. The crude functional equation for the local coefficient of split classical groups over function fields, Theorem~5.14 of [\emph{loc. cit.}], reads
\begin{equation*}
  L^S(s,\tau \times \pi) L^S(2s,\tau,r \circ \rho_m) = \prod_{v \in S} C_{\psi_v}(s,\tau_v \otimes \tilde{\pi}_v,w_0) 
  							    L^S(1-s,\tilde{\tau} \times \tilde{\pi}) L^S(1-2s,\tilde{\tau},r \circ \rho_m),
\end{equation*}
where $r = {\rm Sym}^2$ or $\wedge^2$ depending on the classical group and we use the conventions of \cite{lomelipreprint}. For every $v \in S$ we have that
\begin{equation*}
   C_{\psi_v} = \gamma(s,\tau_v \times \pi_v,\psi_v) \gamma(2s,\tau_v,r \circ \rho_m,\psi_v).
\end{equation*}
The study of exterior and symmetric square $\gamma$-factors begun in \cite{lomeli2009} is completed in \cite{hl2011,lomelipreprint}. We have the functional equation:
\begin{equation*}
   L^S(s,\tau,r \circ \rho_m) = \prod_{v \in S} \gamma(s,\tau_v,r \circ \rho_m, \psi_v) L^S(1-s,\tilde{\tau},r \circ \rho_m).
\end{equation*}
Combining this equation with the crude functional equation of the local coefficient gives property~(vii) for the corresponding $\gamma$-factors:
\begin{equation*}
   L^S(s,\tau \times \pi) = \prod_{v \in S} \gamma(s,\tau_v \times \pi_v,\psi_v) L^S(1-s,\tilde{\tau} \times \tilde{\pi}).
\end{equation*}

Therefore, we can conclude that the existence part of Theorem~\ref{gammathm} holds in the case of $\mathfrak{ls}(p,{\rm GL}_m,{\bf G}_n)$. Notice that no assumption on the characteristic is made for this part of the theorem.

\subsection{The general case} A system of $\gamma$-factors on $\mathfrak{ls}(p,{\rm GL}_m,{\bf  G}_n)$ gives a system of $\gamma$-factors on $\mathfrak{ls}(p,{\bf G}_m,{\rm GL}_n)$ via the relationship
\begin{equation*}
   \gamma(s,\pi \times \tau,\psi) \mathrel{\mathop:} = \gamma(s,\tau \times \pi,\psi),
\end{equation*}
for $(F,\pi,\tau,\psi) \in \mathfrak{ls}(p,{\bf G}_m,{\rm GL}_n)$. We now build upon \S~9 of \cite{lomeli2009}, which is written under the assumption $p \neq 2$. We make this assumption for the rest of this section. Also, we focus on the case of two classical groups ${\bf G}_m$ and ${\bf G}_n$; we let $M$ and $N$ denote the positive integers obtained from $m$ and $n$ via the table on \S~1.3.

First, let $(F,\tau,\pi,\psi) \in \mathfrak{ls}(p,{\bf G}_m,{\bf G}_n)$ be such that $\tau$ is supercuspidal. By Proposition~9.4 of \cite{lomeli2009}, there exists a generic representation ${\rm T}$ of ${\rm GL}_M$ such that
\begin{equation*}
   \gamma(s,\tau \times \rho,\psi) = \gamma(s,{\rm T} \times \rho,\psi),
\end{equation*}
for every supercuspidal representation $\rho$ of ${\rm GL}_r(F)$. The representation ${\rm T}$ is unique due to Th\'eor\`eme~1.1 of \cite{henniart1993}; it is called the local functorial lift of $\tau$.

An irreducible generic discrete series representation $\tau$ of a classical group can be described in terms of its inducing data
\begin{equation} \label{discreteinduced}
   \tau \hookrightarrow {\rm ind}_{P_m}^{G_m} (\delta_1 \otimes \cdots \otimes \delta_d
   									  \otimes \delta_1' \otimes \cdots \otimes \delta_e' \otimes \tau_0),
\end{equation}
where the $\delta_i$'s and the $\delta_j'$'s are essentially square integrable representations of ${\rm GL}_{m_i}(F)$ and $\tau_0$ is an irreducible generic supercuspidal representation of $G_{m_0}$, with ${\bf G}_{m_0}$ is a classical group of the same type as ${\bf G}_m$. Following the results of M\oe glin-Tadi\'c \cite{mt2002}, this is made precise in equation~(9.1) of \cite{lomeli2009}. The case of $m_0 = 0$ is allowed, with appropriate interpretations for the corresponding formulas. If ${\rm T}_0$ is the local functorial lift of $\tau_0$, then the local functorial lift ${\rm T}$ of $\tau$ is the generic constituent of
\begin{equation} \label{discreteimage}
   {\rm ind}_{P_m}^{G_m} (\delta_1 \otimes \cdots \otimes \delta_d \otimes \delta_1' \otimes \cdots \otimes \delta_e' {\rm T}_0 \otimes 
   					  \tilde{\delta_e'} \otimes \cdots \otimes \tilde{\delta_1'} \otimes \tilde{\delta}_d \cdots \otimes \tilde{\delta}_1).
\end{equation}
The representation ${\rm T}$ is a self-dual tempered representation of ${\rm GL}_M(F)$.

Now, an irreducible generic tempered representation of $G_m$ is a constituent of
\begin{equation} \label{temperedinduced}
   {\rm ind}_{P_m}^{G_m} (\delta_1 \otimes \cdots \delta_d \otimes \tau_0),
\end{equation}
where the $\delta_i$'s are discrete series representations of ${\rm GL}_{m_i}(F)$ and $\tau_0$ is a generic discrete series of $G_{m_0}$, where ${\bf G}_{m_0}$ is a classical group of the same kind as ${\bf G}_m$. Then, the local functorial lift ${\rm T}$ of $\tau$ is given by
\begin{equation} \label{temperedimage}
   {\rm ind}_{P_m}^{G_m} (\delta_1 \otimes \cdots \otimes \delta_d \otimes {\rm T}_0 \otimes \tilde{\delta}_d \otimes \tilde{\delta}_1),
\end{equation}
where ${\rm T}_0$ is the local functorial lift of the generic discrete series representation $\tau_0$.

An arbitrary irreducible generic representation $\tau$ of $G_m$ can be described via the work of Mui\'c on the standard module conjecture \cite{muic2001}. Then,
\begin{equation} \label{generalinduced1}
   \tau = {\rm ind}_{P_m}^{G_m} (\tau_1 \nu^{r_1} \otimes \cdots \otimes \tau_{d} \nu^{r_d} \otimes \tau_0).
\end{equation}
Here, each $\tau_i$ is a tempered representation of ${\rm GL}_{m_i}(F)$; $\tau_0$ is a generic tempered representation of $G_{m_0}$, where ${\bf G}_{m_0}$ is a classical group of the same kind as ${\bf G}_m$; and, $\nu = \left| {\rm det}(\cdot) \right|_F$. If ${\bf G}_m = {\rm SO}_{2n}$, and $\tau_0$ is the trivial representation of ${\bf G}_{m_0}(F)$ and $m_d = 1$, the above formula needs the following modification
\begin{equation} \label{generalinduced2}
   \tau = {\rm ind}_{P_m}^{G_m} (\tau_1 \nu^{r_1} \otimes \cdots \otimes \tau_d \nu^{r_d}),
\end{equation}
where we have $0 < \left| r_d \right| < r_{d-1} < \cdots < r_1$. In all other cases, it is given by equation~\eqref{generalinduced1}, where the exponents can be taken so that $0 < r_d < \cdots < r_1$. The local functorial lift ${\rm T}$ of $\tau$ is then given as the unique irreducible quotient of
\begin{equation} \label{generalimage}
   {\rm ind}_{P_M}^{{\rm GL}_M(F)} (\tau_1 \nu^{r_1} \otimes \cdots \otimes \tau_d \nu^{r_d} \otimes {\rm T}_0 \otimes
   							 \tilde{\tau}_d \nu^{-r_d} \otimes \cdots \otimes \tilde{\tau}_1 \nu^{r_1} ),
\end{equation}
where ${\rm T}_0$ is the local functorial lift of $\tau_0$, with appropriate modifications if the induced representation is given by \eqref{generalinduced2}. The local lift has the property that
\begin{equation*}
   \gamma(s,\tau \times \rho,\psi) = \gamma(s,{\rm T} \times \rho,\psi),
\end{equation*}
for any irreducible generic representation $\rho$ of ${\rm GL}_r(F)$.

With a description of the local image of functoriality, we can now obtain a system of extended $\gamma$-factors. Given $(F,\tau,\pi,\psi) \in \mathfrak{ls}(p,{\bf G}_m,{\bf G}_n)$ let $(F,{\rm T},\pi,\psi) \in \mathfrak{ls}(p,{\rm GL}_M,{\bf G}_n)$ be such that ${\rm T}$ is the local functorial lifts of $\tau$. Then, we define
\begin{equation} \label{extendedgamma}
   \gamma(s,\tau \times \pi,\psi) \mathrel{\mathop:}= \gamma(s,{\rm T} \times \pi,\psi).
\end{equation}
We have that $\gamma(s,\tau \times \pi,\psi) = \gamma(s,\pi \times \tau,\psi)$ for every $(F,\tau,\pi,\psi) \in \mathfrak{ls}(p,{\bf G}_m,{\bf G}_n)$. Furthermore, if $(F,{\rm T},\Pi,\psi) \in \mathfrak{ls}(p,{\rm GL}_M,{\rm GL}_N)$, where ${\rm T}$ and $\Pi$ are the local functorial images of $\tau$ and $\pi$, then
\begin{equation*}
   \gamma(s,\tau \times \pi,\psi) = \gamma(s,{\rm T} \times \Pi,\psi).
\end{equation*}

It is now an exercise to show that properties (i) and (ii) are verified; property~(iii) remains as before; and, our definition is indeed compatible with multiplicativity, property (iv). The dependence on $\psi$ can now be obtained from  the corresponding property for $\gamma(s,{\rm T} \times \Pi,\psi)$ (see property~(iv) in the main Theorem of \cite{hlappendix}). And, stability remains as before.

Given $(k,\tau,\pi,\psi) \in \mathcal{LS}(p,{\bf G}_m,{\bf G}_n)$ let ${\rm T}$ and ${\rm \Pi}$ be the functorial lifts of $\tau$ and $\pi$. This is possible via Theorem~9.1 of \cite{lomeli2009}. In fact, the work of Ginzburg, Rallis and Soudry allows us to give a precise description of the image of functoriality \cite{soudry2005}. The global functorial lift ${\rm T}$ of $\tau$ can be expressed as an isobaric sum
\begin{equation} \label{globalisobaric}
   {\rm T} = {\rm T}_1 \boxplus \cdots \boxplus {\rm T}_d,
\end{equation}
where each ${\rm T}_i$, $1 \leq i \leq d$, is a unitary self-dual cuspidal automorphic representation of ${\rm GL}_{M_i}(\mathbb{A}_k)$. Also, $\Pi_i \ncong \Pi_j$ whenever $i \neq j$. Furthermore, if $S$ is a sufficiently large finite set of places of $k$, then
\begin{itemize}
   \item[(i)] $L^S(s,\Pi_i,\wedge^2 \rho_m)$ has a pole at $s=1$, if ${\bf G}_m = {\rm SO}_{2m+1}$;
   \item[(ii)] $L^S(s,\Pi_i,{\rm Sym}^2 \rho_m)$ has a pole at $s=1$, if ${\bf G}_m = {\rm SO}_{2m}$ or ${\rm Sp}_{2m}$.
\end{itemize}
We can similarly express the global functorial lift $\Pi$ of $\pi$ as an isobaric sum.

The functional equation for extended $\gamma$-factors can then be obtained from the above description of the global image and the functional equation for Rankin-Selberg $\gamma$-factors, i.e.,
\begin{align*}
   L^S(s,\tau \times \pi) & = L^S(s,{\rm T} \times \Pi) \\
  & = \prod_{v \in S} \gamma(s,{\rm T}_v \times \Pi_v,\psi_v) L^S(1-s,\tilde{{\rm T}} \times \tilde{\Pi}) \\
  & = \prod_{v \in S} \gamma(s,\tau_v \times \pi_v,\psi_v) L^S(1-s,\tau \times \pi),
\end{align*}
for every $(k,\tau,\pi,\psi) \in \mathcal{LS}(p,{\bf G}_m,{\bf G}_n)$.

\section{Uniqueness}

In the cases involving classical groups, we use a variation of a local-to-global result of Vign\'eras, which follows from the proof of Theorem~2.2 of \cite{v}. We note that, over a global function field, a place at infinity plays the role that archimedean places play over number fields; the notion of an automorphic representation over function fields is independent of the choice of place at infinity. To prove the following proposition, we start with a local field $F$ and take a global field $k$ such that $k_{v_0} \simeq F$ at a place $v_0$ of $k$. We fix two different places $v_1$ and $v_2$ over the same function field $k$. Then one can apply the observation made on p.~469 of [\emph{loc.\,cit.}] to the construction of globally generic cuspidal automorphic representations $\tau$ and $\pi$ from the local representations $\tau_0$ and $\pi_0$. Throughout this section, we impose no restriction on $p$.

\subsection{Proposition} \label{vlocalglobal} \emph{Let $(F,\tau_0,\pi_0,\psi_0) \in \mathfrak{ls}(p,{\bf G}_m,{\bf G}_n)$ be supercuspidal. Then, there exists a $(k,\tau,\pi,\psi) \in \mathcal{LS}(p,{\bf G}_m,{\bf G}_n)$ and a set of places $S = \left\{ v_0, v_1, v_2 \right\}$ of $k$ such that 
\begin{itemize}
   \item[(i)] $k_{v_0} \simeq F$;
   \item[(ii)] $\tau_{v_0} \simeq \tau_0$ and $\pi_{v_0} \simeq \pi_0$;
   \item[(ii)] $\tau_v$ is an unramified principal series for $v \notin \left\{ v_0, v_1 \right\}$;
   \item[(iv)] $\pi_v$ is an unramified principal series for $v \notin \left\{ v_0, v_2 \right\}$. 
\end{itemize}
}

\subsection{Uniqueness for $\mathfrak{ls}(p,{\rm GL}_1,{\bf G}_n)$} \label{GL1classical} Let $\gamma$ be a rule on $\mathfrak{ls}(p,{\rm GL}_1,{\bf G}_n)$ which assigns to every $(F,\chi,\pi,\psi) \in \mathfrak{ls}(p,{\rm GL}_1,{\rm G}_n)$ a rational function on $q_F^{-s}$ satisfying properties (i)-(vii). Using property~(iv), we can reduce to the case when $\pi$ is supercuspidal.

Given $(F,\chi_0,\pi_0,\psi_0) \in \mathfrak{ls}(p,{\rm GL}_1,{\bf G}_n)$ supercuspidal we can lift it to a global $(k,\chi,\pi,\psi) \in \mathcal{LS}(p,{\rm GL}_1,{\bf G}_n)$ where $\pi_v$ is an unramified principal series for $v \notin \left\{ v_0, v_1 \right\}$ as in Proposition~3.1. However, in this situation we can take a character $\chi_{v_0}$ of $k_{v_0}^\times$ which is isomorphic to $\chi_0$ and a highly ramified character $\chi_{v_2}$ of $k_{v_2}^\times$. We then apply the Grundwald-Wang theorem of class field theory \cite{artintate}, in order to lift $\chi_{v_0}$ and $\chi_{v_2}$ to a gr\"ossencharakter $\chi: k^\times \backslash \mathbb{A}_k^\times \rightarrow \mathbb{C}^\times$. From properties (i) and (ii) we have that
\begin{equation*}
   \gamma(s,\chi_0 \times \pi_0,\psi_0) = \gamma(s,\chi_{v_0} \times \pi_{v_0},\psi_{v_0}),
\end{equation*}
and we can assume $\psi_{v_0}$ is obtained from $\psi_0$ using property (v) if necessary.

For every $v \notin \left\{ v_0, v_2 \right\}$, we have that
\begin{equation*}
   {\rm ind}_{B_n}^{G_n}(\mu_{1,v} \otimes \cdots \otimes \mu_{n,v}),
\end{equation*}
where $\mu_{1,v}$, \ldots, $\mu_{n,v}$ are unramified characters. If ${\bf G}_n = {\rm SO}_{2n}$ or ${\rm SO}_{2n+1}$, then
\begin{equation} \label{gl1unramified1}
   \gamma(s,\chi_v \times \pi,\psi) = \prod_{i=1}^n \gamma(s,\chi_v \mu_{i,v},\psi) \gamma(s,\chi_v \mu_{i,v}^{-1}).
\end{equation}
And, if ${\bf G}_n = {\rm Sp}_{2n}$, then
\begin{equation} \label{gl2unramified2}
   \gamma(s,\chi_v \times \pi,\psi) = \gamma(s,\chi_v,\psi) \prod_{i=1}^n \gamma(s,\chi_v \mu_{i,v},\psi) \gamma(s,\chi_v \mu_{i,v}^{-1}).
\end{equation}
The $\gamma$-factors on the right hand side of the previous two equations are abelian $\gamma$-factors of class field theory. Hence, the rule $\gamma$ is uniquely determined at these places.

At the place $v_2$, we let $\xi_{1}$, \ldots $\xi_{n}$, be characters of ${\rm GL}_1(F)$ such that the restriction of $\xi_1 \otimes \cdots \otimes \xi_n$ to the center of $G_n$ agrees with the central character $\omega_{\pi}$ of $\pi$. Let $\tau_{v_2}$ be the generic constituent of
\begin{equation*}
   {\rm ind}_{B_n}^{G_n}(\xi_1 \otimes \cdots \xi_n).
\end{equation*}
Since we have $\chi_{v_2}$ sufficiently ramified we can use property~(vi) to obtain
\begin{equation*}
   \gamma(s,\chi_{v_2} \times \pi_{v_2},\psi_{v_2}) = \gamma(s,\chi_{v_2} \times \tau_{v_2},\psi_{v_2}).
\end{equation*}
Then, using multiplicativity, $\gamma(s,\chi_{v_2} \times \tau_{v_2},\psi_{v_2})$ can be written as a product of abelian $\gamma$-factors similar to equations~\eqref{gl1unramified1} and \eqref{gl2unramified2} above. Any system of $\gamma$-factors satisfying properties (i)-(vi) of the Theorem gives the same result at $v_2$.

Now, let $S$ be a finite set of places of $k$ including $v_0$ and such that $\chi_v$, $\pi_v$, and $\psi_v$ are unramified for $v \notin S$. Then, property (vii) gives
\begin{equation} \label{fegl1case}
   L^S(s,\tau \times \pi) =  \gamma(s,\tau_{v_0} \times \pi_{v_0},\psi_{v_0}) 
   					\prod_{v \in S - \left\{ v_0 \right\}} \gamma(s,\tau_v \times \pi_v,\psi_v) \, 
					L^S(1-s,\tilde{\tau}_v \times \tilde{\pi}_v).
\end{equation}
Since $\gamma$-factors are determined for every $v \notin S - \left\{ v_0 \right\}$ by the above discussion, we can conclude that $\gamma(s,\tau_{v_0} \times \pi_{v_0},\psi_{v_0})$ is completely determined by equation~\eqref{fegl1case}.

\subsection{Uniqueness in general} Let $\gamma$ be a rule on $\mathfrak{ls}(p)$ which assigns to every quadruple $(F,\tau,\pi,\psi) \in \mathfrak{ls}(p)$ a rational function on $q_F^{-s}$ satisfying properties (i)-(vii). Using property~(iv), we can reduce to the supercuspidal case.

Take a fixed supercuspidal $(F,\tau_0,\pi_0,\psi_0) \in \mathfrak{ls}(p,{\bf G}_m,{\bf G}_n)$ and lift it to a global $(k,\tau,\pi,\psi) \in \mathcal{LS}(p,{\bf G}_m,{\bf G}_n)$ via Proposition~\ref{vlocalglobal}, properties~(i), (ii) and (v). Let ${\bf B}_m$ (resp. ${\bf B}_m$) be the Borel subgroup of ${\bf G}_m$ (resp. ${\bf G}_n$) of upper triangular matrices. For every $v \notin \left\{ v_0, v_1 \right\}$, let $\chi_{1,v}$, \ldots, $\chi_{m,v}$ be unramified characters of ${\rm GL}_1(k_v)$ such that $\tau_v$ occurs as a subrepresentation of the unitarily induced representation
\begin{equation*}
   {\rm ind}_{B_m}^{G_m}(\chi_{1,v} \otimes \cdots \otimes \chi_{m,v}).
\end{equation*}
For every $v \notin \left\{ v_0, v_2 \right\}$, let $\mu_{1,v}$, \dots, $\mu_{n,v}$ be unramified characters of ${\rm GL}_1(k_v)$ such that $\pi_v$ occurs as a subrepresentation of the unitarily induced representation
\begin{equation*}
   {\rm ind}_{B_n}^{G_n}(\mu_{1,v} \otimes \cdots \otimes \mu_{n,v}).
\end{equation*}
Take $v \notin S$, then both $\tau_v$ and $\pi_v$ are unramified. We can in then use properties (iii) and (iv) to show that:
\begin{itemize}
   \item[(a)] If both ${\bf G}_m$ and ${\bf G}_n$ are classical groups, then
   \begin{align*}
      \gamma(s,\tau_v \times \pi_v,\psi_v)
      		= \prod_{i=1}^d &\gamma(s,\chi_{i,v} \mu_{0,v},\psi_v) \gamma(s,\chi^{-1}_{i,v} \mu_{0,v},\psi_v)
      		      \prod_{j=1}^e \gamma(s,\chi_{0,v} \mu_{j,v},\psi_v) \gamma(s,\chi_{0,v} \mu^{-1}_{j,v},\psi_v) \\
		\times \prod_{1 \leq h \leq d, 1 \leq l \leq e} & \gamma(s,\chi_{h,v} \mu_{l,v},\psi_v) \gamma(s,\chi_{h,v}\mu^{-1}_{l,v},\psi_v) \gamma(s,\chi^{-1}_{h,v} \mu_{l,v},\psi_v) \gamma(s,\chi^{-1}_{h,v} \mu^{-1}_{l,v},\psi_v),
   \end{align*}
   where we take $\gamma(s,\chi_{0,v} \mu_{j,v},\psi_v)$ and $\gamma(s,\chi_{0,v} \mu^{-1}_{j,v},\psi_v)$ (resp. $\gamma(s,\chi_{i,v} \mu_{0,v},\psi_v)$ and $\gamma(s,\chi^{-1}_{i,v} \mu_{0,v},\psi_v)$) to be trivial if ${\bf G}_m$ (resp. ${\bf G}_n$) is a special orthogonal group; and, we take $\chi_0 = 1$ (resp. $\mu_0 = 1$) if ${\bf G}_m$ (resp. ${\bf G}_n$) is symplectic. 
   \item[(b)] If ${\bf G}_m = {\rm GL}_m$ and ${\bf G}_n$ is a classical group, then
   \begin{equation*}
      \gamma(s,\tau_v \times \pi_v,\psi_v)
      		= \prod_{i=0}^d \gamma(s,\chi_{i,v} \mu_{0,v},\psi_v) 
      		\times \prod_{i=1}^d \prod_{j=1}^e \gamma(s,\chi_{i,v} \mu_{j,v},\psi_v) \gamma(s,\chi_{i,v} \mu_{j,v}^{-1},\psi_v),
   \end{equation*}
   where we take $\gamma(s,\chi_{j,v} \mu_{0,v},\psi_v)$ to be trivial if ${\bf G}_n$ is a special orthogonal group; and, we take $\mu_0 = 1$ if ${\bf G}_n$ is a symplectic group.
   \item[(c)] If ${\bf G}_m = {\rm GL}_m$ and ${\bf G}_n = {\rm GL}_n$, then
   \begin{equation*}
      \gamma(s,\tau_v \times \pi_v,\psi_v) = \prod_{i,j} \gamma(s,\chi_{i,v} \mu_{j,v},\psi_v).
   \end{equation*}
   \end{itemize}
Any system of $\gamma$-factors satisfying properties (i)-(iv) gives $\gamma(s,\tau_v \times \pi_v,\psi_v)$, $v \notin S$, as a product of abelian $\gamma$-factors of Tate's thesis as above. Hence, it is uniquely determined at these places.

The remaining possibility, at either $v_1$ or $v_2$, is that one representation is unramified while the other remains arbitrary. For concreteness, assume $\tau_{v_1}$ and $\pi_{v_2}$ are unramified while $\tau_{v_2}$ and $\pi_{v_1}$ remain arbitrary. Then $\tau_{v_1}$ and $\pi_{v_2}$ are constituents of representations via unitary parabolic induction from a product of unramified characters as before. The multiplicativity property of $\gamma$-factors gives $\gamma(s,\tau_{v_1} \times \pi_{v_1},\psi_{v_1})$ as a product of $\gamma$-factors of the form $\gamma(s,\chi_{v_1} \times \pi_{v_1},\psi_{v_1})$, where $\chi_{v_1}$ is a character of ${\rm GL}_1(k_{v_1})$. Similarly, multiplicativity gives an expression for $\gamma(s,\tau_{v_2} \times \pi_{v_2},\psi_{v_2})$ as a product of $\gamma$-factors of the form $\gamma(s,\tau_{v_2} \times \mu_{v_2})$, where $\mu_{v_2}$ is a character of ${\rm GL}_1(k_{v_2})$. In these cases, $\gamma$-factors are uniquely determined as shown in \S~\ref{GL1classical}, where property~(vi) is used.

At places where $\psi_v$ may be ramified, we can use property~(v) to compute $\gamma$-factors with respect to an unramified character. The functional equation for $\gamma$-factors gives,
\begin{equation*}
   L^{S'}(s,\tau \times \pi) =	\gamma(s,\tau_{v_0} \times \pi_{v_0},\psi_{v_0}) 
   					\prod_{v \in S' - \left\{ v_0 \right\}} \gamma(s,\tau_{v} \times \pi_v,\psi_v) \, 
					L^{S'}(1-s,\tilde{\tau}_v \times \tilde{\pi}_v).
\end{equation*}
Since partial $L$-functions are uniquely determined, and we have shown that $\gamma$-factors are uniquely determined at places other than $v_0$, we can conclude that $\gamma(s,\tau_{v_0} \times \pi_{v_0},\psi_{v_0})$ is uniquely determined.

\section{Extended $L$-functions and root numbers}

We now turn towards the defining properties of $L$-functions and $\varepsilon$-factors. We assume that $p \neq 2$, which is necessary to study the case $\mathfrak{ls}(p,{\bf G}_m,{\bf G}_n)$, when ${\bf G}_m$ and ${\bf G}_n$ are both classical groups.

\subsection{Axioms for a local system of $L$-functions and root numbers} With a system of $\gamma$-factors on $\mathfrak{ls}(p)$ satisfying properties (i)-(vii), we can proceed to define rational functions $L(s,\tau \times \pi)$ and monomials $\varepsilon(s,\tau \times \pi,\psi)$ on the variable $q_F^{-s}$ for every $(F,\tau,\pi,\psi) \in \mathfrak{ls}(p)$.

\begin{itemize}
   \item[(viii)] (Tempered $L$-functions). \emph{For $(F,\tau,\psi,\psi) \in \mathfrak{ls}(p)$ tempered, let $P_{\tau \times \pi}(t)$ be the unique polynomial with $P_{\tau \times \pi}(0) = 1$ and such that $P_{\tau \times \pi}(q_F^{-s})$ is the numerator of $\gamma(s,\tau \times \pi,\psi)$. Then
   \begin{equation*}
      L(s,\tau \times \pi) = P_{\tau \times \pi}(q^{-s})^{-1}.
   \end{equation*}
   }
   \item[(ix)] (Holomorphy of tempered $L$-functions). \emph{Let $(F,\tau,\pi,\psi) \in \mathfrak{ls}(p)$ be tempered. Then, $L(s,\tau \times \pi)$ is holomorphic and non-zero for $\Re(s) > 0$.}
    \item[(x)] (Tempered $\varepsilon$-factors). \emph{Let $(F,\tau,\pi,\psi) \in \mathfrak{ls}(p)$ be tempered, then
   \begin{equation*}
      \varepsilon(s,\tau \times \pi,\psi) = \gamma(s,\tau \times \pi,\psi) \dfrac{L(s,\tau \times \pi)}{L(1-s,\tilde{\tau} \times \tilde{\pi})}.
   \end{equation*}
   }
   \item[(xi)] (Twists by unramified characters). \emph{Let $(F,\tau,\pi,\psi) \in \mathfrak{ls}(p,{\rm GL}_m,{\bf G}_n)$, then
   \begin{align*}
      L(s+s_0,\tau \times \pi) & = L(s,\left| \det(\cdot) \right|_F^{s_0} \tau \times \pi) \\
      \varepsilon(s+s_0,\tau \times \pi,\psi) & = \varepsilon(s,\left| \det(\cdot) \right|_F^{s_0} \tau \times \pi,\psi).
   \end{align*}
   }
   \item[(xii)] (Multiplicativity). \emph{Let $(F,\tau_i,\pi_j,\psi) \in \mathfrak{ls}(p,{\rm GL}_{m_i},{\rm GL}_{n_j})$, $1 \leq i \leq d$, $1 \leq j \leq e$, be quasi-tempered, and let $(F,\tau_0,\pi_0,\psi) \in \mathfrak{ls}(p,{\bf G}_{m_0},{\bf G}_{n_0})$ be tempered. Let $m = \sum m_i$ and $n = \sum n_j$. Let ${\bf G}_m$ and ${\bf G}_n$ be of the same type as ${\bf G}_{m_0}$ and ${\bf G}_{n_0}$. Let ${\bf P}_m$ and ${\bf P}_n$ be parabolic subgroups of ${\bf G}_m$ and ${\bf G}_n$ with Levi subgroups ${\bf M}_m \simeq \prod_{i=1}^d {\rm GL}_{m_i} \times {\bf G}_{m_0}$ and ${\bf M}_n \simeq \prod_{j=1}^e {\rm GL}_{n_j} \times {\bf G}_{n_0}$. Suppose that $(F,\tau,\pi,\psi) \in \mathfrak{ls}(p)$ is such that $\tau$ is the generic constituent of
   \begin{equation*}
      {\rm ind}_{P_m}^{G_m} (\tau_1 \otimes \cdots \otimes \tau_d \otimes \tau_0),
   \end{equation*}
   and $\pi$ is the generic constituent of
   \begin{equation*}
      {\rm ind}_{P_n}^{G_n} (\pi_1 \otimes \cdots \otimes \pi_e \otimes \pi_0).
   \end{equation*}
   When $m_0=0$ (resp. $n_0=0)$ we make the following conventions: take $\tau_0$ (resp. $\pi_0$) to be the trivial character of ${\rm GL}_1(F)$ if ${\bf G}_m$ (resp. ${\bf G}_n$) is a symplectic group; in all other cases, we interpret local factors involving $\tau_0$ (resp. $\pi_0$) to be trivial. 
   \begin{itemize}
   \item[(xii.a)] If both ${\bf G}_m$ and ${\bf G}_n$ are classical groups, then
\begin{align*}
   L&(s,\tau \times \pi) = L(s,\tau_0 \times \pi_0) \\
   				& \times \prod_{i=1}^d L(s,\tau_i \times \pi_0) L(s,\tilde{\tau}_i \times \pi_0)
   				   \prod_{j=1}^e L(s,\tau_0 \times \pi_j) L(s,\tau_0 \times \tilde{\pi}_j) \\
				& \times \prod_{i=1}^d \prod_{j=1}^e L(s,\tau_i \times \pi_j) L(s,\tilde{\tau}_i \times \pi_j) 
				   L(s,\tau_i \times \tilde{\pi}_j) L(s,\tilde{\tau}_i \times \tilde{\pi}_j),
\end{align*}
and local root numbers satisfy
\begin{align*}
   \varepsilon&(s,\tau \times \pi,\psi) = \varepsilon(s,\tau_0 \times \pi_0,\psi) \\
   	& \times \prod_{i=1}^d \varepsilon(s,\tau_i \times \pi_0,\psi) \varepsilon(s,\tilde{\tau}_i \times \pi_0,\psi) 
	   \prod_{j=1}^e \varepsilon(s,\tau_0 \times \pi_j,\psi) \varepsilon(s,\tau_0 \times \tilde{\pi}_j,\psi) \\
   	& \times \prod_{i=1}^d \prod_{j=1}^e \varepsilon(s,\tau_i \times \pi_j,\psi) \varepsilon(s,\tilde{\tau}_i \times \pi_j,\psi)
	   \varepsilon(s,\tau_i \times \tilde{\pi}_j,\psi) \varepsilon(s,\tilde{\tau}_i \times \tilde{\pi}_j,\psi).
\end{align*}
   \item[(xii.b)] If ${\bf G}_m = {\rm GL}_m$ and ${\bf G}_n$ is a classical group, then
   \begin{equation*}
      L(s,\tau \times \pi) = \prod_{i=0}^d L(s,\tau_i \times \pi_0) 
      				      \prod_{i=1}^d \prod_{j=1}^e L(s,\tau_i \times \pi_j) L(s,\tau_i \times \tilde{\pi}_j),
   \end{equation*}
   and
   \begin{equation*}
      \varepsilon(s,\tau \times \pi,\psi) = \prod_{i=0}^d \varepsilon(s,\tau_i \times \pi_0,\psi)
      		\prod_{i=1}^d \prod_{j=1}^e \varepsilon(s,\tau_i \times \pi_j,\psi) \varepsilon(s,\tau_i \times \tilde{\pi}_j,\psi).
   \end{equation*}
   \item[(xii.c)] If ${\bf G}_m = {\rm GL}_m$ and ${\bf G}_n = {\rm GL}_n$, then
   \begin{equation*}
      L(s,\tau \times \pi) = \prod_{i=0}^d \prod_{j=0}^e L(s,\tau_i \times \pi_i),
   \end{equation*}
   and
   \begin{equation*}
      \varepsilon(s,\tau \times \pi,\psi) = \prod_{i=0}^d \prod_{j=0}^e \varepsilon(s,\tau_i \times \pi_j,\psi).
   \end{equation*}
   \end{itemize}
   }
\end{itemize}

\subsection{Additional properties of $\gamma$-factors} The following properties are satisfied by any system of $\gamma$-factors on $\mathfrak{ls}(p)$ satisfying properties (i)-(vii).

\begin{itemize}
   \item[(xiii)] (Local functional equation). \emph{Let $(F,\tau,\pi,\psi) \in \mathfrak{ls}(p)$, then
   \begin{equation*}
      \gamma(s,\tau \times \pi,\psi)\gamma(1-s,\tilde{\tau} \times \tilde{\pi},\overline{\psi}) = 1.
   \end{equation*}
   }
   \item[(xiv)] (Twists by unramified characters). \emph{Let $(F,\tau,\pi,\psi) \in \mathfrak{ls}(p,{\rm GL}_m,{\bf G}_n)$, then
   \begin{equation*}
      \gamma(s+s_0,\tau \times \pi,\psi) = \gamma(s,\left| \det(\cdot) \right|_F^{s_0} \tau \times \pi,\psi).
   \end{equation*}
   }
   \item[(xv)] (Commutativity). \emph{Let $(F,\tau,\pi,\psi) \in \mathfrak{ls}(p)$, then
   \begin{equation*}
      \gamma(s,\tau \times \pi,\psi) = \gamma(s,\pi \times \tau,\psi).
   \end{equation*}
   }
\end{itemize}

We now give a proof of property~(xiii) following the proof of uniqueness given in \S~3. Notice that it is a property of abelian $\gamma$-factors: if $(F,\chi,\mu,\psi) \in \mathfrak{ls}(p,{\rm GL}_1,{\rm GL}_1)$, then
\begin{equation*}
   \gamma(s,\chi \mu ,\psi) \gamma(1-s,\chi^{-1}\mu^{-1},\overline{\psi}) = 1,
\end{equation*}
see for example equation~(1.3) of \cite{lomelipreprint}. First, we prove the local functional equation for the case of $\mathfrak{ls}(p,{\rm GL}_1,{\bf G}_n)$. We can reduce to the supercuspidal case via multiplicativity. Lift a supercuspidal $(F,\chi_0,\pi_0,\psi_0) \in \mathfrak{ls}(p,{\rm GL}_1,{\bf G}_n)$ to $(k,\chi,\pi,\psi) \in \mathcal{LS}(p,{\rm GL}_1,{\bf G}_n)$ as in \S~3.2. Notice that the method of proof gives an expression for every $\gamma(s,\chi_v \times \pi_v,\psi_v)$, $v \notin \left\{ v_0 \right\}$, as a product of abelian $\gamma$-factors. Then, the local functional equation at $v_0$ follows from the global functional equation applied twice.

Now, in the proof of uniqueness for the general case of \S~3.3, let $(k,\tau,\pi,\psi) \in \mathcal{LS}(p,{\bf G}_m,{\bf G}_n)$ be the global quadruple obtained from the supercuspidal quadruple $(F,\tau_0,\pi_0,\psi_0) \in \mathfrak{ls}(p,{\bf G}_m,{\bf G}_n)$. The method of proof gives an expression for every $\gamma(s,\tau_v \times \pi_v,\psi_v)$, $v \notin \left\{ v_0 \right\}$, in terms of $\gamma$-factors for $\mathfrak{ls}(p,{\rm GL}_1,{\bf G}_n)$, $\mathfrak{ls}(p,{\bf G}_m,{\rm GL}_1)$ or $\mathfrak{ls}(p,{\rm GL}_1,{\rm GL}_1)$. In all of these cases, the local functional equation holds. Hence, it follows at $v_0$ by applying the global functional equation twice.

We leave the proofs of properties (xiv) and (xv) as exercises.

\subsection{Theorem} \emph{It $p \neq 2$, there {\bf exists} a system of local factors on $\mathfrak{ls}(p)$ satisfying properties (i)-(xii). Any system of local factors on $\mathfrak{ls}(p)$ satisfying properties (i)-(xii) is {\bf uniquely} determined.}

\subsection*{Proof} We have already established the existence and uniqueness part of the theorem concerning properties (i)-(vii). We now construct local $L$-functions and $\varepsilon$-factors.

We first treat the tempered case, where property~(viii) is taken as the definition of a local $L$-function. Notice that the multiplicativity property of $\gamma$-factors gives the multiplicativity property of local $L$-functions in the tempered case.

We now prove property~(ix) for the new case of two classical groups ${\bf G}_m$ and ${\bf G}_n$. Let $(F,\tau,\pi,\psi) \in \mathfrak{ls}(p,{\bf G}_m,{\bf G}_n)$ be tempered. The representation $\tau$ is a constituent of
\begin{equation*}
   {\rm ind}_{P_m}^{G_m} (\delta_1 \otimes \cdots \otimes \delta_d \otimes \tau_0),
\end{equation*}
as in equation~\eqref{temperedinduced} with $\delta_i$, $i = 1$, \ldots $d$, generic discrete series representations of ${\rm GL}_{m_i}(F)$, and $\tau_0$ a generic discrete series representation of $G_{m_0}$. Similarly, $\pi$ is a constituent of
\begin{equation*}
   {\rm ind}_{P_n}^{G_n} (\rho_1 \otimes \cdots \otimes \rho_e \otimes \pi_0),
\end{equation*} 
with $\rho_j$, $j = 1$, \ldots, $e$, generic discrete series representations of ${\rm GL}_{m_j}(F)$, and $\pi_0$ a generic discrete series representation of $G_{n_0}$.

Let ${\rm T}_0$ and $\Pi_0$ be the local functorial lifts of $\tau_0$ and $\pi_0$ given by equation~\eqref{discreteimage}. Notice that they are self-dual tempered representations of ${\rm GL}_{M_0}(F)$ and ${\rm GL}_{N_0}(F)$. The local functorial lift ${\rm T}$ of $\tau$ is given by
\begin{equation*}
   {\rm ind}_{P_m}^{G_m} (\delta_1 \otimes \cdots \otimes \delta_d \otimes {\rm T}_0 \otimes \tilde{\delta}_1 \otimes \cdots \otimes \tilde{\delta}_1)
\end{equation*}
and the lift $\Pi$ of $\pi$ is given by
\begin{equation*}
   {\rm ind}_{P_n}^{G_n} (\rho_1 \otimes \cdots \otimes \rho_e \otimes \Pi_0 \otimes \tilde{\rho}_e \otimes \cdots \otimes \tilde{\rho}_1).
\end{equation*}
Then
\begin{align*}
   L(s,\tau \times \pi) = & L(s,{\rm T}_0 \times \Pi_0) \prod_{i=1}^d L(s,\delta_i \times \Pi_0) L(s,\tilde{\delta}_i \times \Pi_0)
   				      \prod_{j=1}^e L(s,{\rm T}_0 \times \rho_j) L(s,{\rm T}_0 \times \tilde{\rho}_j) \\
   				   & \times \prod_{i=1}^d \prod_{j=1}^e L(s,\delta_i \times \rho_j) L(s,\tilde{\delta}_i \times \rho_j)
   				   L(s,\delta_i \times \tilde{\rho}_j) L(s,\tilde{\delta}_i \times \tilde{\rho}_j).
\end{align*}
Each $L$-function on the right hand side is a Rankin-Selberg $L$-function for products of representations of general linear groups, known to be holomorphic for $\Re(s) > 0$. Hence, the extended $L$-function $L(s,\tau \times \pi)$ is holomorphic for $\Re(s) > 0$. Thus, our local $L$-functions satisfy property~(ix) in the tempered case.

Next, for tempered $(F,\tau,\pi,\psi) \in \mathfrak{ls}(p)$, property~(x) is taken as the definition of root numbers. Then $\varepsilon(s,\tau \times \pi,\psi)$ is a monomial in $q_F^{-s}$. For this, we use the local functional equation of $\gamma$-factors, property~(xiii), together with property~(ix) to ensure that no cancellations occur on the strip $0 < \Re(s) < 1$.

Having defined local $L$-functions and root numbers for tempered representations, they are then defined on $\mathfrak{ls}(p)$ with the aid of Langlands classification. In the case of ${\bf G}_m = {\rm GL}_m$ and ${\bf G}_n = {\rm GL}_n$, we include a treatment of $L$-functions and $\varepsilon$-factors in \cite{lomelipreprint} in a self contained manner within the $\mathcal{LS}$-method. The definition given by equations (7.6) of [\emph{loc. cit.}] is in accordance with the definition of Rankin-Selberg local factors \cite{jpss1983}. We will now define extended $L$-functions and root numbers on $\mathfrak{ls}(p,{\bf G}_m,{\bf G}_n)$, when both ${\bf G}_m$ and ${\bf G}_n$ are classical groups. Obtaining explicit defining relations for local $L$-functions and $\varepsilon$-factors on $\mathfrak{ls}(p,{\rm GL}_m,{\bf G}_n)$, with ${\bf G}_n$ a classical group, is left as an exercise.

Assume that both ${\bf G}_m$ and ${\bf G}_n$ are classical groups. Let $\tau$ be given by equation~\eqref{generalinduced1}
\begin{equation*}
   {\rm ind}_{P_m}^{G_m} (\tau_1 \nu^{r_1} \otimes \cdots \otimes \tau_{d} \nu^{r_d} \otimes \tau_0),
\end{equation*}
where each $\tau_i$ is a tempered representation of ${\rm GL}_{m_i}(F)$ and $\tau_0$ is a generic tempered representation of $G_{m_0}$. With the appropriate modifications if $\tau$ is given by equation~\eqref{generalinduced2}. Similarly, $\pi$ is given by
\begin{equation*}
   {\rm ind}_{P_n}^{G_n} (\pi_1 \nu^{t_1} \otimes \cdots \otimes \pi_{e} \nu^{t_e} \otimes \pi_0),
\end{equation*}
where each $\pi_i$ is a tempered representation of ${\rm GL}_{n_i}(F)$ and $\pi_0$ is a generic tempered representation of $G_{n_0}$. With appropriate modifications if $\pi$ is given by equation~\eqref{generalinduced2}.

Let $(F,\tau,\pi,\psi) \in \mathfrak{ls}(p)$, then we define local $L$-functions by
\begin{align*}
   L&(s,\tau \times \pi) = L(s,\tau_0 \times \pi_0) \\
   				& \times \prod_{i=1}^d L(s+r_i,\tau_i \times \pi_0) L(s-r_i,\tilde{\tau}_i \times \pi_0)
   				   \prod_{j=1}^e L(s+t_j,\tau_0 \times \pi_j) L(s-t_j,\tau_0 \times \tilde{\pi}_j) \\
				& \times \prod_{i=1}^d \prod_{j=1}^e L(s + r_i + t_j,\tau_i \times \pi_j) L(s-r_i+t_j,\tilde{\tau}_i \times \pi_j) \\
				& \times L(s+r_i-t_j,\tau_i \times \tilde{\pi}_j) L(s-r_i-r_j,\tilde{\tau}_i \times \tilde{\pi}_j),
\end{align*}
and local root numbers by
\begin{align*}
   \varepsilon&(s,\tau \times \pi,\psi) = \varepsilon(s,\tau_0 \times \pi_0,\psi) \\
   	& \times \prod_{i=1}^d \varepsilon(s+r_i,\tau_i \times \pi_0,\psi) \varepsilon(s-r_i,\tilde{\tau}_i \times \pi_0,\psi) 
	   \prod_{j=1}^e \varepsilon(s+t_j,\tau_0 \times \pi_j,\psi) \varepsilon(s-t_j,\tau_0 \times \tilde{\pi}_j,\psi) \\
   	& \times \prod_{i=1}^d \prod_{j=1}^e \varepsilon(s+r_i+t_j,\tau_i \times \pi_j,\psi) \varepsilon(s-r_i+t_j,\tilde{\tau}_i \times \pi_j,\psi) \\
	& \times \varepsilon(s+r_i-t_i,\tau_i \times \tilde{\pi}_j,\psi) \varepsilon(s-r_i-t_j,\tilde{\tau}_i \times \tilde{\pi}_j,\psi).
\end{align*}

It is now possible to show that our construction is compatible with properties (xi) and (xii). Notice that the definition of local $L$-functions and root numbers is based on an extended system of $\gamma$-factors, and only uses special cases of properties (viii)-(xii). Hence, any system of extended $\gamma$-factors, local $L$-functions and root numbers satisfying properties (i)-(xii) is uniquely determined.

\subsection{Theorem} \emph{Assume that $p \neq 2$. For every $(k,\tau,\pi,\psi) \in \mathcal{LS}(p)$, let
\begin{equation} \label{globalLepsilon}
   L(s,\tau \times \pi) = \prod_{v} L(s,\tau_v \times \pi_v), \quad 
   \varepsilon(s,\tau \times \pi) = \prod_{v} \varepsilon(s,\tau_v \times \pi_v,\psi_v).
\end{equation}
Automorphic $L$-functions satisfy the following properties:}
\begin{itemize}
   \item[(i)] (Rationality). \emph{$L(s,\tau \times \pi)$ converges on a right half plane and has a meromorphic continuation to a rational function on $q^{-s}$.}
   \item[(ii)] (Functional equation). $L(s,\tau \times \pi) = \varepsilon(s,\tau \times \pi) L(1-s,\tilde{\tau} \times \tilde{\pi})$.
   \item[(iii)] (Riemann Hypothesis) \emph{The zeros of $L(s,\tau \times \pi)$ are contained in the line $\Re(s) = 1/2$}.
\end{itemize}

\subsection*{Proof} First, given $(F,\tau,\pi,\psi) \in \mathfrak{ls}(p,{\rm GL}_m,{\bf G}_n)$, we know that partial $L$-functions $L^S(s,\tau \times \pi)$ converge on a right half plane. The rationality argument of Harder for Eisenstein series \cite{harder1974}, allows us to give an automorphic forms proof that $L^S(s,\tau \times \pi)$ is a rational function on $q^{-s}$, Corollary~6.6 of \cite{lomeli2009}. Now, notice that each local $L$-function in the product
\begin{equation*}
   \prod_{v \in S} L(s,\tau_v \times \pi_v)
\end{equation*}
is a rational function on $q_v^{-s} = q^{-{\rm deg}(v)s}$. Hence, property~(i) follows for completed automorphic $L$-functions in this case. Also, the definition of local $L$-functions and $\varepsilon$-factors at the places $v \in S$ can be incorporated into the functional equation satisfied by $\gamma$-factors in order to obtain property~(ii) for global $L$-functions and $\varepsilon$-factors on $\mathcal{LS}(p,{\rm GL}_m,{\bf G}_n)$.

Next, we treat the case of $\mathcal{LS}(p,{\bf G}_m,{\bf G}_n)$, with both ${\bf G}_m$ and ${\bf G}_n$ classical groups. Let $(k,\tau,\pi,\psi) \in \mathcal{LS}(p,{\bf G}_m,{\bf G}_n)$. The functorial lift of Theorem~9.1 of \cite{lomeli2009} is compatible with the local functorial lift at every place $v$ of $k$. The construction of the local lift is reviewed in \S~2.4. Equation~\eqref{globalisobaric} enables us to write the global lifts ${\rm T}$ of $\tau$ and $\Pi$ of $\pi$ as isobaric sums
\begin{equation} \label{globallifts}
   {\rm T} = {\rm T}_1 \boxplus \cdots \boxplus {\rm T}_d \text{  and  } \Pi = \Pi_1 \boxplus \cdots \boxplus \Pi_e.
\end{equation}
We have that
\begin{equation*}
   L^S(s,\tau \times \pi) = L^S(s,{\rm T} \times \Pi),
\end{equation*}
which converge on a right half plane and have a meromorphic continuation to a rational function on $q^{-s}$. At the remaining places, extended $L$-functions $L(s,\tau_v \times \pi_v)$ are rational on $q_v^{-s}$. Hence, the completed automorphic $L$-function $L(s,\tau \times \pi)$ satisfies property~(i). The way local factors are defined can be coupled with the functional equation satisfied by extended $\gamma$-factors in order to establish property~(ii) for automorphic $L$-functions on $\mathcal{LS}(p,{\bf G}_m,{\bf G}_n)$.

Finally, let $(k,\tau,\pi,\psi) \in \mathcal{LS}(p,{\bf G}_m,{\bf G}_n)$. If both ${\bf G}_m$ and ${\bf G}_n$ are classical groups, take ${\rm T}$ and $\Pi$ to be the global functorial lifts of $\tau$ and $\pi$ of equation~\eqref{globallifts}. If ${\bf G}_m = {\rm GL}_m$, we just take ${\rm T} = \tau$. Then
\begin{equation*}
   L(s,\tau \times \pi) = L(s,{\rm T} \times \Pi) = \prod_{i,j} L(s,{\rm T}_i \times \Pi_j),
\end{equation*}
where $(k,{\rm T}_i,\Pi_j,\psi) \in \mathcal{LS}(p,{\rm GL}_{m_i},{\rm GL}_{n_j})$, for $1 \leq i \leq d$, $1 \leq j \leq e$. Thanks to the work of Lafforgue on the global Langlands conjecture for ${\rm GL}_N$ over function fields, each Rankin-Selberg $L$-function $L(s,{\rm T}_i \times \Pi_j)$ satisfies the Riemann Hypothesis (see Th\'eor\`eme~VI.10(ii) of \cite{laff2002}). Hence, so does $L(s,\tau \times \pi)$. We conclude that automorphic $L$-functions on $\mathcal{LS}(p)$ satisfy the Riemann Hypothesis.

\end{document}